\documentclass[12pt]{article}
\def\be {\begin{equation}}
\def\ee {\end{equation}}
\def\ba {\begin{eqnarray}}
\def\ea {\end{eqnarray}}

\def\bi {\begin{itemize}}
\def\ei {\end{itemize}}
\begin{document}
\def\bea{\begin{eqnarray}}
\def\eea{\end{eqnarray}}
\title{\textbf{Solution for a non linear Schrodinger equation via Hopf-Cole transformation} }
\author{  \textbf{H. R. Rezazadeh},\textbf{ Tirdad Soulati}}

\maketitle

\begin{abstract}
In this short letter we show that the  one dimensional non linear
Schr\"{o}dinger equation(NLS)can be solved by a Hopf-Cole
transformation which converts it to the Burgers equation in
turbulence.
 \end{abstract} \maketitle

\section{General remarks}
NLS often outbreaks in various fields of the applied physics.Many
authors both mathematician and physicists investigated it's
analytical properties\cite{group}.For a NLS, one can not obtains a
compact solution even when we restrict ourselves only to a simple
dimensional reduction of the independent variables to 2,i.e only
$x,t$.Recently a  NLS proposed which new non linear term is
gravitational potential that arises  only due to the entropic force
acting on a particle with the aid of the entropy defined from
quantum mechanics\cite{Joakim}.

As a simple direct evidence of the seem of a NLS,suppose that
a free  particle with mass $m$ moves in flat spacetime.
We know that in some scenarios of the early universe, Planck
constant spotted as an effective momentum dependence function ,via
Generalized Uncertainty Principle (GUP)\cite{GUP}
 and not constant regardless on some
further developments in non commutative geometry. The idea that we
introduce here can be staminate to the mass and regard it as an
effective function dependence on the momentum of
position\cite{mass}.The simple presumption may be a linear
dependence to the momentum in the $\hbar$ or $m$.Everybody can treat
with it as an effective Schr\"{o}dinger equation for a test
particle(free in rest) that a linear  stimulus with unknown quantum
essence disturb it's mass or further some GUP corrections have been
replaced the usual Planck constant with a momentum dependence
effective one. One applicable is the following NLS,
\begin{eqnarray}
i\hbar\dot{\psi}=-\frac{\hbar^{2}}{2m}\psi''-i\hbar
l^{\frac{1}{2}}c\psi'\psi
\end{eqnarray}
with $\psi(x,t),\dot{f}=\partial_{t}f,f'=\partial_{x}f$.It is
salutary  for reader if we  alight dimensions of this new term. We
take the light speed $c=1$(Geometrical units),and for accuracy of
the usual dimension of the wave function in one dimensional problems
we introduced a unit length scale $l=1$(can be regarded as the
 inverse square of the norm of the wave function
 i.e.$l^{-1}=\sqrt{||\psi||_{C^{2}}}$.We read (1) as the entropic
 efficacy of the Holographic screen of the test particle in the
 neighborhood screen in the new explanation  of the Gravity as an
 entropic force exerted on the particle near a holographic
 screen(equip potential surfaces) in an emergence model of the
 spacetimes \cite{Verlinde}.\\
\subsection{Exact solution}
 Anyway the equation (1) is a
simple NLSE,and by definition
\begin{eqnarray}
\varepsilon=\frac{i\hbar}{2m}
\end{eqnarray}

 converts to the following famous equation in non linear
PDE, Burgers' eq\cite{Burgers'},
\begin{eqnarray}
\dot{\psi}=\varepsilon \psi''-\psi\psi'
\end{eqnarray}
Burgers' equation is a fundamental partial differential equation
from fluid mechanics. It occurs in various areas of applied
mathematics, such as modeling of gas dynamics and traffic flow. It
is named for Johannes Martinus Burgers'. Now consider the planar 2
dimensional vector field,
\begin{eqnarray}
V=(-u,0.5 \psi^{2}-\varepsilon \psi')
\end{eqnarray}
Obviously is Curl free, i.e. a conservative vector field. Thus
according to the familiar theorems, there exists a scalar potential
function which obeys the next conditions,
\begin{eqnarray}
\phi'=-\psi,\dot{\phi}=0.5\psi^{2}-\varepsilon \psi'
\end{eqnarray}
The new scalar function $\phi$ solves the equation
\begin{eqnarray}
\dot{\phi}=0.5\phi'^{2}+\varepsilon \phi''
\end{eqnarray}
Letting $\phi=2\varepsilon \log(\eta)$,then we are left with
\begin{eqnarray}
\dot{\eta}-\varepsilon \eta''=0
\end{eqnarray}
The one dimensional Heat equation  or free particle Schr\"{o}dinger
wave equation and from (5) we obtain
\begin{eqnarray}
\psi=-2\varepsilon\frac{\eta'}{\eta}
\end{eqnarray}
Which is the Hopf-Cole transformation\cite{HC}.The Burgers' equation
then  been linearized by the Cole-Hopf transformation.Some solutions
of the Burgers equation has been investigated by Majid and
Ranasinghe \cite{Majid} .The general solution for (3) with Initial
Conditions (IC) $\psi_{0}(x)=\psi(x,0)$ is
\begin{eqnarray}
\psi(x,t)=\frac{\int_{-\infty}^{\infty}\frac{x-y}{t}\psi_{0}(y)K(x,y;t)dy}{\int_{-\infty}^{\infty}K(x,y;t)dy}
\end{eqnarray}
Where in it the heat kernel or propagator for a free particle in
usual can be written as,
\begin{eqnarray}
K(x,y;t)=exp(\frac{im(x-y)^{2}}{2t\hbar})
\end{eqnarray}
WE can related this propagator to the partition function of a
statistical system by replacing the temperature with a Wick rotated
time scale.
\subsection{Localized particle}
Now for a particle which is localized at origin at time $t=0$ we can
substitute the initial non renormalizable wave function

\begin{eqnarray}
\psi_{0}(x)=N\delta(x)
\end{eqnarray}
We know that the wave function for a localized particle is not
renormalizable.In some texts the authors set $N=1$.But we exempt N
from this, regard it as a free constant. and if we assume that for
initial data we must have
\begin{eqnarray}
_{|x|\rightarrow\infty}\frac{1}{x^{2}}\int^{x}_{a}\psi_{0}(y)dy\rightarrow0
\end{eqnarray}
(a is arbitrary and never affect the value of the wave
function,choosing a=1)and after simple integration we have
\begin{eqnarray}
\psi(x,t)=\sqrt{\frac{2i\hbar}{m\pi
t}}\frac{exp(\frac{imx^{2}}{2\hbar
t})}{2(exp(-imN/\hbar)-1)^{-1}+\frac{\sqrt{\pi}}{2}(1-erf(\frac{x\sqrt{m}}{\sqrt{2i\hbar
t}})}
\end{eqnarray}
where $ erf(x)\equiv\int^{x}_{0}e^{-z^{2}}dz$ is the integral of the
Gaussian distribution .
As \cite{Joakim}  we suggest that the probability density
$|\psi(x,t) |$ is in fact related to a partition function $Z$ for
different possible states,in the locution of the statistical
mechanics.Thus we derived the partition function $Z$ for different
possible quantum states in the NLS evolutionary scheme.
\section{Summary}
NLS equations arisen in different theoretical and applied
problems.One of the most class of the solvable models due to the
Burgers'.In this work we look on the non linear terms in the
Burgers' equation as the entropic corrections to the free particle
Schr\"{o}dinger equation in the frame of the new model proposed by
Verlinde\cite{Verlinde}.We show that under influence of this non
linear term,the time evolution  of a localized particle  is so
different from the common quantum mechanics.We related the norm of
the wave function to the partition function of the system and also
explicitly we   calculate the wave function for it.

\end{document}